\documentclass[12pt]{article}
\usepackage{amssymb}

\usepackage{wrapfig}
\usepackage[dvips]{graphicx}
\usepackage{graphicx}
\usepackage{amsmath}
\usepackage{amsmath}

\newtheorem{theorem}{Theorem}

\newtheorem{lemma}[theorem]{Lemma}

\input{tcilatex}

\begin{document}

\title{Slowing down and reflection of waves in truncated periodic media.}
\author{S. Molchanov and \ B. Vainberg \footnote{Corresponding author,
E-mail address: brvainbe@uncc.edu  (B. Vainberg). Both authors
were partially supported by the NSF grant DMS-0405927. } \and
Dept. of Mathematics, University of North Carolina at Charlotte,
\and Charlotte, NC 28223, USA}
\date{}
\maketitle

\begin{abstract}
We consider wave propagation and scattering governed by one dimensional Schr%
\"{o}dinger operators with truncated periodic potentials. The propagation of
wave packets with narrow frequency supports is studied. The goal is to
describe potentials for which the group velocity (for periodic problem) is
small and the transmission coefficient for the truncated potential is not
too small, i.e. to find media where a slowing down of the wave packets
coexists with a transparency.


\end{abstract}

\textbf{I. Introduction}

This paper is designed to provide a mathematical foundation for studying the
processes of slowing down of the wave packets in periodic media. The problem
of creating compact and efficient optical delay devices has been very
actively discussed in literature for the last several years, but it has not
yet been solved technologically. Such devices can find many applications,
for example in synchronizing the work of very fast optical elements and much
slower electronics. Several different mechanisms to slow down the
propagation of wave packets have been considered (see \cite{y1}, \cite{y2},
\cite{y3}, \cite{hb}, \cite{hb2}, \cite{f}, \cite{mv}, \cite{mv1}). Most of
them are based on a use of periodic media with a flat dispersion relation $%
\omega =\omega (k)$ in some region. For such media, the derivative $%
k^{\prime }(\omega )$ is very big, and the group velocity $V_{g}(\omega )=%
\frac{L}{k^{\prime }(\omega )}$ is small for a specific interval of
frequencies $|\omega -\omega _{0}|<\varepsilon $. Here $L$ is the period.

Let $H$ be a periodic 1-D Schr\"{o}dinger operator on $L^{2}(R):$

\begin{equation}
H\psi =-\frac{d^{2}\psi }{dx^{2}}+q(x)\psi ,\text{ \ \ \ }q(x+L)=q(x).
\label{ha}
\end{equation}
Let $\psi =\psi (x,\omega )$ be a Bloch solution of the equation $H\psi
=\omega ^{2}\psi ,$ i.e. $\psi (x+L,\omega )=\psi (x,\omega )e^{ik(\omega )}.
$ The function $k=k(\omega )$ defines the dispersion relation, and $e^{\pm
ik(\omega )}$ are the eigenvalues of the monodromy operator $M_{\omega }$.
The function $\psi (x,\omega )e^{-i\omega t}$ is a solution of the
nonstationary wave equation $u_{tt}=-Hu$, it corresponds to a wave which
propagates along the $x$-axis with the (phase) velocity $V(\omega )=\omega
L/k(\omega ).$ Consider a wave packet
\begin{equation}
u=\underset{\left| \omega -\omega _{0}\right| <\varepsilon }{\int }\psi
(x,\omega )\alpha (\omega -\omega _{0})e^{-i\omega t}d\omega   \label{pa}
\end{equation}
with angular frequencies from a small interval around $\omega _{0}.$ The
group velocity for the packet (\ref{pa}) is defined as $V_{g}=L/k^{\prime
}(\omega _{0}).$ One can easily show that
\begin{equation}
u\sim \gamma (x)\widetilde{\alpha }(\frac{x}{V_{g}}-t)e^{i(k(\omega _{0})%
\frac{x}{L}-\omega _{0}t)}  \label{asf}
\end{equation}
where $\widetilde{\alpha }$ is the Fourier transform of the function $\alpha
(\omega ),$ $\gamma $ is a periodic function, i.e. $u$ has the form of a
plane wave which propagates with the phase velocity $V(\omega _{0})$ and
which is modulated by another wave $\widetilde{\alpha }$ whose speed of
propagation is $V_{g}.$ Asymptotic formula (\ref{asf}) is valid when $%
\varepsilon $ and $\varepsilon ^{2}\left| x\right| $ are small. The main
physical idea behind the optical delay devices is to find an optical
material for which the group velocity $V_{g}$ is very small on some interval
around a specific ''singular'' frequency $\omega _{0}$ and then to use the
wave $\widetilde{\alpha }$ (rather than the last factor in the right hand
side of (\ref{asf})) for the transmission of information.

It is well known that the spectrum of the periodic operator $H$ has a
band-gap structure with the bands defined by the equation $|$Trace $%
M_{\omega }|\leq 2.$ The function $k(\omega )$ changes
monotonically on any band, the changes are equal to $\pm \pi$, and
therefore, $k^{\prime }$ is large and $V_{g}$ is small on any
narrow band. This suggests solving the problem of slowing down of
the wave packets by creating a periodic optical material with
narrow spectral bands. A group of physicists from Caltech (A.
Yariv and others \cite{y1}, \cite{y2}, \cite{y3}) proposed to use
a special coupled resonators optical waveguide (CROW) as an
optical delay device. It can be realized as a set of identical
periodic cavities (arranged in a linear sequence) in a dielectric
medium. If there are no interactions between resonators, then the
spectrum of the system consists of discrete eigenvalues (of the
problem in an individual cavity) of infinite multiplicity. If a
weak and periodic interaction between neighboring resonators is
imposed, then such eigenvalues will generate narrow optical bands.
P. Kuchment and A. Kunyansky \cite{kuku} provided spectral
analysis of very general high contrast periodic media which
yielded many examples of high contrast media with narrow spectral
bands (and therefore, with the flat dispersion relations and small
group velocity).

Hence, one can expect that periodic media with small group
velocity (for example, due to the existence of a narrow spectral
band) may serve as optical delay devices. However, let us note
that all the arguments above are related to infinite periodic
media, and a real device has finite size. The main goal of this
paper is to answer the following question: in what sense may the
problem in an infinite periodic medium serve as an approximation
for the processes of wave propagation through a long finite slab
of the periodic medium. This question becomes particularly
important due to the following circumstance: one needs to know if
the incident pulses with the frequencies supported in a small
neighborhood of a singular frequency $\omega _{0}$\ (where the
dispersion is flat) will enter the slab or, perhaps, the majority
of their energy will be reflected.

In this paper we study a Hamiltonian $H_{N}$ of the form
\begin{equation}
H_{N}\psi =-\frac{d^{2}\psi }{dx^{2}}+q_{N}(x)\psi ,  \label{han}
\end{equation}
where the potential $q_{N}$ coincides with a periodic real valued potential $%
q(x)$ when $x\in \lbrack 0,NL]$, and $q_{N}=0$ outside the interval $%
[0,NL],$ and $N$ is very big. The goal is to find a potential $q$
such that for some fixed interval $\Delta $ on the frequency axis
the group velocity for operator (\ref{ha}) is small and the
transmission coefficient for operator (\ref{han}) is not too
small. We shall see that in order to achieve this goal we need to
consider a family of operators with potentials $q(x)=Av(x),$
$A\rightarrow \infty .$ The parameter $A$ has to be large in order
to make the group velocity small. At the same time the interval
$\Delta $ shrinks as $A\rightarrow \infty $, and therefore, we
need to assume that $N\rightarrow \infty $ (a physical motivation
for that will be given). We shall see that there are two scenarios
with similar results: $\Delta $ is a narrow band for the operator
(\ref{ha}) or it is a neighborhood of a degenerated band edge. The
dispersion relation is very flat in the first case, and it is not
flat in the second case. In both cases, the length of $\Delta $
and the group velocity have order $O(A^{-1})$. The transmission
coefficient has the same order on the main part of the interval
$\Delta $, but there are $N$ very narrow transparency zones on
$\Delta .$ We shall also justify the fact that in spite of these
transparency zones, only a small part of the energy (of order
$O(A^{-2})$) of the incident pulse will be transmitted through the
media, and the majority of the energy will be reflected.

We shall not discuss other known models of slowdown of wave
packets: inflection point (A. Figotin, I. Vitebskii, \cite{f}),
necklace model (S. Molchanov, B. Vainberg \cite{mv}), SCISSOR (J.
Heebner, R. Boyd, \cite{hb},\cite{hb2}) mainly because the
corresponding mathematical tools are different from the classical
1-D scalar Schr\"{o}dinger operator: it is symplectic systems in
\cite{f}, quantum graphs in \cite{mv}, etc.

The paper has the following structure.

In the second section, we develop the scattering theory for the
operator ( \ref{han}) combining the ideas of the Floquet-Bloch
theory and the 1-D scattering theory for operators with fast
decaying potentials. In particular, exact formulas for the
reflection and transmission coefficients will be obtained through
the monodromy operator and $N.$ Section II contains also a
transition to semi-infinite periodic materials ($N\rightarrow
\infty $). The relation between a slowing down and the reflection
is discussed in the last section.

\textbf{II. Scattering by a finite slab of a periodic medium.}

a). \textit{General Schr\"{o}dinger operators and operators with
periodic or compactly supported potentials. }We shall first recall
some well known facts concerning the spectral problem for the
general Schr\"{o}dinger operators on $L^{2}\left( R\right) :$
\begin{equation}
H\psi =-\psi ^{\prime \prime }+q\left( x\right) \psi =\omega ^{2}\psi ,
\label{s0}
\end{equation}
where the potential $q$ (which can be a distribution)\ is real
valued and bounded from below in an appropriate sense, say, there
is a constant $c$ such that
\begin{equation}
\int_{x}^{x+1}\min (q\left( x\right) ,0)dx>c\text{ \ \ for all \ }x\in R.
\label{par}
\end{equation}
The last inequality provides the uniqueness of the self adjoint extension in
$L^{2}\left( R\right) $ of the operator $H$ defined originally on the space
of infinitely smooth, compactly supported functions $\psi .$

\bigskip Let $T_{\omega }\left( 0,x\right) $ be the Pr\"{u}ffer transfer
matrix (propagator) for operator $H$
\begin{equation}
T_{\omega }\left( 0,x\right) =\left[
\begin{array}{ll}
\psi _{1}\left( \omega ,x\right) & \omega \psi _{2}\left( \omega ,x\right)
\\
\frac{\psi _{1}^{\prime }\left( \omega ,x\right) }{\omega } & \psi
_{2}^{\prime }\left( \omega ,x\right)
\end{array}
\right]  \label{s1}
\end{equation}
where $\psi _{1,2}\left( \omega ,x\right) $ are the solutions of the
equation $H\psi =\omega ^{2}\psi $ with initial data
\begin{equation}
\psi _{1}\left( \omega ,0\right) =1,\text{ \ }\psi _{1}^{\prime
}\left( \omega ,0\right) =0;\text{ \ \ \ \ }\psi _{2}\left( \omega
,0\right) =0, \text{ \ }\psi _{2}^{\prime }\left( \omega ,0\right)
=1,  \label{s2}
\end{equation}
(i.e., simply, $T_{\omega }\left( 0,0\right) $ is the identity
matrix). For any solution $\psi $ of the equation $H\psi =\omega
^{2}\psi ,$ matrix $ T_{\omega }\left( 0,x\right) $ maps the
Pr\"{u}ffer Cauchy data of $\psi $ at $x=0$ into the Pr\"{u}ffer
Cauchy data of $\psi $ at point $x:$
\begin{equation*}
T_{\omega }\left( 0,x\right) :\left[
\begin{array}{l}
\psi \left( 0\right) \\
\frac{\psi ^{\prime }\left( 0\right) }{\omega }
\end{array}
\right] \rightarrow \left[
\begin{array}{l}
\psi \left( x\right) \\
\frac{\psi ^{\prime }\left( x\right) }{\omega }
\end{array}
\right] .
\end{equation*}

Assumption (\ref{par}) implies that the equation $H\psi =\omega ^{2}\psi $
with Im$\omega >0$ has exactly one solution $\psi ^{+}$ in $L^{2}\left(
R_{+}\right) $ normalized by the condition $\psi ^{+}\left( \omega ,0\right)
=1,$ and it has exactly one solution $\psi ^{-}\in L^{2}\left( R_{-}\right) $
normalized by the same condition. Here $R_{\pm }$ are the semiaxes $%
x\gtrless 0.$ \ Obviously, $\psi ^{\pm }$ can be represented as linear
combinations of $\psi _{1}$ and $\psi _{2},$\ and from the normalization of $%
\psi ^{\pm }$ it follows that there exist functions $m^{\pm }=m^{\pm }\left(
\omega \right) $ such that
\begin{equation*}
\psi ^{\pm }=\psi _{1}+m^{\pm }\left( \omega \right) \psi _{2},\text{ \ \ Im}%
\omega >0.
\end{equation*}
In other words, $m^{\pm }$ are functions such that
\begin{equation*}
\psi _{1}+m^{+}\left( \omega \right) \psi _{2}\in L^{2}\left( R_{+}\right) ,%
\text{ \ }\psi _{1}+m^{-}\left( \omega \right) \psi _{2}\in L^{2}\left(
R_{-}\right) ,\text{ \ \ \ Im}\omega >0.
\end{equation*}
Functions $m^{\pm }=m^{\pm }\left( \omega \right) $ are called Weyl's
functions, and $\psi ^{\pm }$ are called Weyl's solutions. \

Let the potential $q=q(x)$ be periodic: $q(x+L)=q(x).$ Consider the
propagator through one period (monodromy matrix):
\begin{equation*}
M_{\omega }=T_{\omega }\left( 0,L\right) =\left[
\begin{array}{ll}
\alpha & \beta \\
\gamma & \delta
\end{array}
\right] \left( \omega \right)
\end{equation*}
Put $F\left( \omega \right) =\dfrac{1}{2}$Tr$M_{\omega }=\left( \psi
_{1}+\psi _{2}^{\prime }\right) \left( \omega ,L\right) =\left( \alpha
+\delta \right) \left( \omega \right) $. Both $M_{\omega }$ and $F\left(
\omega \right) $ are entire functions of $\omega $. Since $\det M_{\omega
}=1,$ the eigenvalues $\mu ^{\pm }\left( \omega \right) $ of $M_{\omega }$
are the roots of the following characteristic equation:
\begin{equation}
\mu ^{2}+2F\left( \omega \right) \mu +1=0.  \label{s3}
\end{equation}

If Im$\omega >0$ then one can select roots $\mu ^{\pm }\left( \omega \right)
$ of (\ref{s3}) in such a way that $\mu ^{\pm }\left( \omega \right) =e^{\pm
ik\left( \omega \right) }$, where $k\left( \omega \right) $ is analytic and
\begin{equation}
\text{Im}k\left( \omega \right) >0\text{ \ when Im}\omega >0,  \label{imk}
\end{equation}
i.e.,
\begin{equation}
\left| \mu ^{+}\left( \omega \right) \right| <1,\text{ \ }\left| \mu
^{-}\left( \omega \right) \right| >1,\text{ \ Im}\omega >0.  \label{aaa}
\end{equation}
The roots $\mu ^{\pm }\left( \omega \right) $ for real $\omega \geq 0$ are
defined by continuity in the upper half plane:
\begin{equation*}
\mu ^{\pm }\left( \omega \right) =\mu ^{\pm }\left( \omega +i0\right) ,\text{
\ \ }\omega \in \left[ 0,\infty \right] .
\end{equation*}
Since the trace of $M_{\omega }$ is equal to the sum of the
eigenvalues $ e^{\pm ik(\omega )},$
\begin{equation}
\cos k(\omega )=F\left( \omega \right) =(\alpha +\delta )/2.  \label{cos}
\end{equation}

If the potential $q$ is not negative then the spectrum of $H$
belongs to the positive part of the energy axis $\lambda =\omega
^{2},$ and the spectrum has a band-gap structure. If $q_{-}=\min
(q,0)\neq 0$ then $H$ may also have spectrum (a finite number of
gaps, in particular) on the negative part of the energy axis
$\lambda =-\omega ^{2}<0.$ Wave processes are associated only with
the positive part of the spectrum, and we shall exclude the
negative part of the spectrum from the future analysis. Thus, the
bands on the frequency axis $\omega $ will always be related to
the positive energies.

For real $\omega \geq 0,$ the inequality $\left| F\left( \omega
\right) \right| \leq 1$ defines the spectral bands (zones)
$b_{n}=\left[ \omega _{2n-1},\omega _{2n}\right] ,$ $n=1,2,...,$\
on the frequency axis $\omega = \sqrt{\lambda }$. The function
$k(\omega )$ is real valued when $\omega $ belongs to a band. The
roots $\mu ^{\pm }\left( \omega \right) $ are complex adjoint
there, and $|\mu ^{\pm }|=1.$ The spectrum of $H$ (on $
L^{2}\left( R\right) $) on the frequency axis is
$\bigcup\limits_{n=1}^{ \infty }b_{n}.$ The complimentary open
set, given by $\left| F\left( \omega \right) \right| >1,$
corresponds to spectral gaps. On gaps, the roots $\mu ^{\pm
}\left( \omega \right) $ are real and (\ref{aaa}) holds. Figure 1
presents a typical graph of $F\left( \omega \right) .$ A point
$\omega _{j} $ which belongs to the boundary of a band and the
boundary of a gap is called a non-degenerate band edge. If it
belongs to the boundary of two different bands, it is called a
degenerate band edge.
\begin{figure}[tbp]
\label{neck1}
\par
\begin{center}
\includegraphics[width=120mm, angle=0]{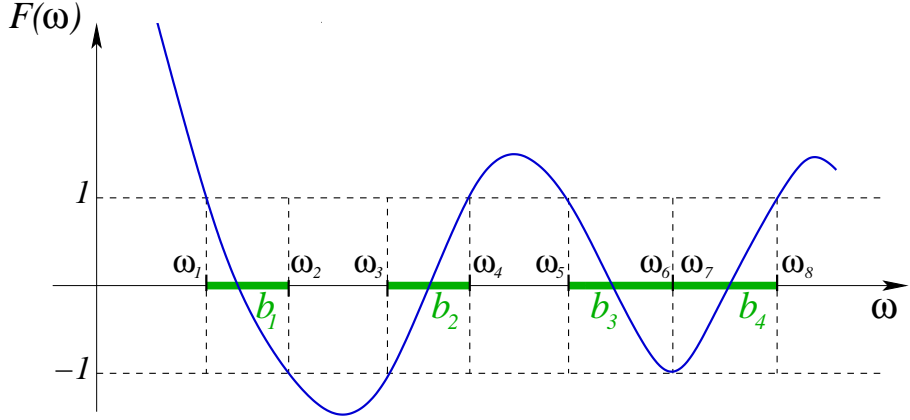}
\end{center}
\caption{A typical graph of $F\left( \protect\omega \right) .$ Segments $%
b_{1}=\left[ \protect\omega _{1},\protect\omega _{2}\right] ,$ $b_{2}=\left[
\protect\omega _{3},\protect\omega _{3}\right] ,$ $b_{3}=\left[ \protect%
\omega _{5},\protect\omega _{6}\right] ,$ $b_{4}=\left[ \protect\omega _{7},%
\protect\omega _{8}\right] $ are bands, the point $\protect\omega _{6}=%
\protect\omega _{7}$ is a degenerate band edge.}
\end{figure}

\begin{lemma}
\label{la}If $\omega =\omega _{0}$ is a non-degenerate band edge,
then $ F^{\prime }(\omega _{0})\neq 0.$ If $\omega =\omega _{0}$
is a degenerate band edge, then $F^{\prime }(\omega _{0})=0,$
$F^{\prime \prime }(\omega _{0})\neq 0.$ Both eigenvalues of the
monodromy matrix $M_{\omega }$ at any band edge are equal to $1$
or both are equal to $-1.$ The matrix $M_{\omega } $ at the
non-degenerate band edges has an off-diagonal element in its
Jordan form, and $M_{\omega }=\pm I $ at the degenerate band
edges, where $I$ is the identity matrix.
\end{lemma}

We shall provide only a sketch of a proof of this Lemma, and we
refer to \cite{ls} for details. Let us recall that bands
are defined by the condition $|F\left( \omega \right) |\leq 1,$ and that $%
F\left( \omega \right) =\pm 1$ at any band edge $\omega =\omega _{0}$. If $%
F^{\prime }(\omega _{0})\neq 0$ then $|F\left( \omega \right) |-1$
changes sign at $\omega =\omega _{0}$, and $\omega =\omega _{0}$
is a
non-degenerate band edge. Let $\omega =\omega _{0}$ be a band edge and $%
F^{\prime }(\omega _{0})=0.$ Using the equation in variations, one
can show that in this case $F^{\prime \prime }(\omega _{0})\neq 0$
and $F(\omega _{0})F^{\prime \prime }(\omega _{0})<0.$ Hence,
$|F\left( \omega \right) |\leq 1$ in a neighborhood of such a band
edge, i.e. $\omega =\omega _{0}$ is a degenerate band edge.

The eigenvalues $\mu ^{\pm }$ of $M_{\omega }$ \ are equal to
$e^{\pm ik\left( \omega \right) }$, where $\cos k\left( \omega
\right) =F\left( \omega \right) =\pm 1$ at a band edge $\omega
=\omega _{0}.$ Hence, both eigenvalues $\mu ^{\pm }$ at $\omega
=\omega _{0}$ are equal to $1$ or both are equal to $-1$ (and
$F(\omega _{0})=1$ or $F(\omega _{0})=-1,$ respectively)$.$ Thus,
$\lambda =$ $\omega _{0}^{2}$ \ is an eigenvalue of
the periodic or, respectively, anti-periodic problem for the operator $%
H $. If the band edge is non-degenerate, then this eigenvalue is
simple, and the Jordan form of $M_{\omega }$ is non-diagonal
(since the eigenspace is one dimensional). If the band edge is
degenerate, then the
corresponding eigenspace is two dimensional, and the Jordan form of $%
M_{\omega }$ is diagonal, and therefore, $M_{\omega }=\pm $.

The proof of Lemma \ref{la} is complete.

We normalize the eigenvectors $h^{\pm }\left( \omega \right) $ of $M_{\omega
}$ by choosing the first coordinate of $h^{\pm }\left( \omega \right) $ to
be equal to one: $h^{\pm }\left( \omega \right) =\left[
\begin{array}{l}
1 \\
m^{\pm }\left( \omega \right)
\end{array}
\right] $. The second coordinates of the vectors $h^{\pm }\left( \omega
\right) $ coincide with Weyl's functions defined above. In fact, if $\psi
^{\pm }$ are solutions of the equation $H\psi =\omega ^{2}\psi $ with the
Pr\"{u}ffer Cauchy data given by the eigenvector $h^{\pm },$ then
\begin{equation}
\psi ^{\pm }\left( \omega ,x+L\right) =e^{\pm ik\left( \omega \right) }\psi
^{\pm }\left( \omega ,x\right) ,  \label{s5}
\end{equation}
and (\ref{aaa}) implies that $\psi ^{\pm }\in $ $L^{2}\left(
R_{\pm }\right) $ when Im$\omega >0.$ From here it follows that
$\psi ^{\pm }$ coincide with Weyl's solutions introduced for
general Hamiltonians $H$, and that the second coordinates of the vectors $%
h^{\pm }\left( \omega \right) $ are Weyl's functions.

Since
\begin{equation}
\left[
\begin{array}{ll}
\alpha -e^{\pm ik\left( \omega \right) } & \text{ \ \ \ \ }\beta \\
\text{ \ \ \ \ }\gamma & \delta -e^{\pm ik\left( \omega \right) }
\end{array}
\right] \left[
\begin{array}{l}
1 \\
m^{\pm }
\end{array}
\right] =0,  \label{s4}
\end{equation}
the following two representations are valid for Weyl's functions:
\begin{equation}
m^{\pm }\left( \omega \right) =\dfrac{e^{\pm ik\left( \omega \right)
}-\alpha \left( \omega \right) }{\beta \left( \omega \right) }=\frac{\gamma
\left( \omega \right) }{e^{\pm ik\left( \omega \right) }-\delta \left(
\omega \right) }.  \label{s4a}
\end{equation}

For real $\omega \in \bigcup\limits_{n=1}^{\infty }b_{n},$ the Weil's
solutions $\psi ^{\pm }\left( \omega ,x\right) $ are the usual Bloch almost
periodic solutions. The relationship between the space and time frequencies $%
k=k\left( \omega \right) ,$ $\omega \in \bigcup\limits_{n=1}^{\infty }b_{n},$
is called the dispersion law.

Consider now a Hamiltonian $H$ with a compactly supported
potential $q$. \ Let supp$q\subseteq \lbrack 0,A].$ A function
$\psi $ is called the solution of the scattering problem for the
operator $H$ if $H\psi =\omega ^{2}\psi $ and
\begin{equation*}
\psi =e^{i\omega x}+re^{-i\omega x},\text{ \ \ }x<0;\text{ \ \ \ \ }\psi
=te^{i\omega x},\text{ \ \ }x>A,
\end{equation*}
with some $r=r(\omega )$ and $t=t(\omega ).$ The numbers $r$ and
$t$ are called the reflection and transmission coefficients,
respectively.

For any matrix $M=[m_{i,j}],$ let $||M||$ be the Hilbert-Schmidt norm of $M:
$%
\begin{equation*}
||M||^{2}=\underset{i,j}{\sum }m_{i,j}^{2}.
\end{equation*}

\begin{theorem}
Let the support of $q\ $belong to $[0,A].$ Then the following relations hold
for the reflection and transmission coefficients:
\begin{equation}
|r(\omega )|^{2}=\frac{||T_{\omega }(0,A)||^{2}-2}{||T_{\omega }(0,A)||^{2}+2%
},\text{ \ \ \ }|t(\omega )|^{2}=\frac{4}{||T_{\omega }(0,A)||^{2}+2}.
\label{p1}
\end{equation}
\end{theorem}

\textbf{Proof.} The Pr\"{u}ffer Cauchy data for the scattering solution at $%
x=0$ and $x=A$ are
\begin{equation*}
\left[
\begin{array}{c}
1+r \\
i\left( 1-r\right)
\end{array}
\right] ,\text{ \ \ }\left[
\begin{array}{c}
t \\
it
\end{array}
\right] e^{i\omega A}\text{,}
\end{equation*}
respectively. This implies
\begin{equation*}
\left[
\begin{array}{c}
t \\
it
\end{array}
\right] e^{i\omega A}=T_{\omega }\left( 0,A\right) \left[
\begin{array}{c}
1+r \\
i\left( 1-r\right)
\end{array}
\right] .
\end{equation*}

We equate the ratios of the second and the first coordinates of
the vectors above. Using the notation
\begin{equation*}
T_{\omega }\left( 0,A\right) =\left[
\begin{array}{ll}
a & b \\
c & d
\end{array}
\right] \left( \omega \right)
\end{equation*}
we arrive at
\begin{equation*}
\frac{c(1+r)+di(1-r)}{a(1+r)+bi(1-r)}=i.
\end{equation*}
Solving for $r$ we obtain
\begin{equation}
r=\frac{d-a-i(c+d)}{a+d+i(c-d)}.  \label{r1}
\end{equation}
Thus,
\begin{equation}
|r|^{2}=\frac{(d-a)^{2}+(c+d)^{2}}{(a+d)^{2}+(c-d)^{2}}=\frac{||T_{\omega
}\left( 0,A\right) ||^{2}-2ad+2cd}{||T_{\omega }\left( 0,A\right)
||^{2}+2ad-2cd}.  \label{p2}
\end{equation}
Since $\det T_{\omega }\left( 0,A\right) =1,$ (\ref{p2}) implies the first
of relations (\ref{p1}). The second relation follows from the first and the
energy conservation law: $|r|^{2}+|t|^{2}=1.$

The proof of the Theorem is complete.

b) \textit{Scattering by truncated periodic potentials. }The main
result of this section concerns the Hamiltonian $H_{N}$ with the
truncated periodic potential $q_{N}:$
\begin{equation*}
H_{N}\psi =-\frac{d^{2}\psi }{dx^{2}}+q_{N}(x)\psi ,\text{ \ \ }%
q_{N}(x)=q(x)I_{[0,NL]},
\end{equation*}
which appears when the propagation of waves through a finite slab
of a periodic medium is studied. We shall also consider the
limiting case $N=\infty :$%
\begin{equation*}
H_{\infty }\psi =-\frac{d^{2}\psi }{dx^{2}}+q_{\infty }(x)\psi ,\text{ \ \ }%
q_{\infty }(x)=q(x)I_{[0,\infty ]},
\end{equation*}
which corresponds to the case of such a long slab that it can be
considered as half infinite.

We shall denote the reflection and transmission coefficients for
the operator $H_{N}$ (with a compactly supported potential
$q_{N})$ by $r_{N}$ and $t_{N}$, respectively. In the case of the
operator $H_{\infty },$ the solution $\psi $ of the scattering
problem is defined as the solution of the equation $H_{\infty
}\psi =\omega ^{2}\psi $ which has the form
\begin{equation}
\psi =e^{i\omega x}+re^{-i\omega x},\text{ \ \ }x<0;\text{ \ \ \ \ }\psi
=c\psi ^{+}(\omega ,x),\text{ \ \ }x>0,  \label{ji}
\end{equation}
with some $r=r(\omega ),$ $c=c(\omega ).$

\begin{theorem}
\label{t1}1). The transfer matrix over $N$ periods $M_{\omega }^{N}$ $%
=T_{\omega }\left( 0,LN\right) $ has the form
\begin{equation}
M_{\omega }^{N}=\left[
\begin{array}{ll}
\alpha _{N} & \beta _{N} \\
\gamma _{N} & \delta _{N}
\end{array}
\right] =\frac{\sin Nk\left( \omega \right) }{\sin k\left( \omega \right) }%
M_{\omega }-\frac{\sin \left( N-1\right) k\left( \omega \right)
}{\sin k\left( \omega \right) }I,  \label{s6}
\end{equation}
where $k=k\left( \omega \right) $ is the dispersion relation. The elements
of $M_{\omega }^{N}$ satisfy the relations
\begin{eqnarray}
\alpha _{N}-\delta _{N} &=&\frac{\sin Nk\left( \omega \right) }{\sin k\left(
\omega \right) }\left( \alpha -\delta \right) ,\text{ \ \ \ }\beta _{N}=%
\frac{\sin Nk\left( \omega \right) }{\sin k\left( \omega \right) }\beta ,
\notag \\
\text{ \ \ }\gamma _{N} &=&\frac{\sin Nk\left( \omega \right) }{\sin k\left(
\omega \right) }\gamma ,\text{ \ \ \ }\alpha _{N}+\delta _{N}=2\cos
Nk(\omega ).  \label{s7}
\end{eqnarray}
The Hilbert-Schmidt norm of $M_{\omega }^{N}$ is equal to
\begin{equation}
||M_{\omega }^{N}||^{2}=(||M_{\omega }||^{2}-2)\frac{\sin ^{2}Nk\left(
\omega \right) }{\sin ^{2}k\left( \omega \right) }+2  \label{ss}
\end{equation}

2). The reflection coefficients have the forms
\begin{eqnarray}
r_{N}\left( \omega \right) &=&-\frac{\left( \alpha -\delta \right) +i\left(
\beta +\gamma \right) }{2\sin k\left( \omega \right) \frac{\cos Nk\left(
\omega \right) }{\sin Nk\left( \omega \right) }+i\left( \gamma -\beta
\right) },\text{ \ \ }  \label{rnr} \\
\text{\ }r\left( \omega \right) &=&\frac{\beta +\gamma -i\left( \alpha
-\delta \right) }{2\sin k\left( \omega \right) +\left( \beta -\gamma \right)
}.  \label{rnr1}
\end{eqnarray}

3). The following formula is valid for $|t_{N}\left( \omega \right) |:$%
\begin{equation}
|t_{N}\left( \omega \right) |^{2}=\frac{4}{(||M_{\omega }||^{2}-2)\frac{\sin
^{2}Nk\left( \omega \right) }{\sin ^{2}k\left( \omega \right) }+4}.
\label{tnt}
\end{equation}

4). The reflection coefficients $r_{N}\left( \omega \right) $ and $r\left(
\omega \right) $ are analytic in the upper half plane $C^{+}=\{\omega :$Im$%
k>0\}$ and continuous in $\overline{C^{+}},$ and $r_{N}\left(
\omega \right) \rightarrow r\left( \omega \right) $ when
$N\rightarrow \infty $ and $\omega \in C^{+}.$ When $\omega $ is
real, $r_{N}\left( \omega \right) $ converges to $r\left( \omega
\right) $ in the weak sense:
\begin{equation*}
\int_{-\infty }^{\infty }r_{N}\left( \omega \right) \varphi (\omega )d\omega
\rightarrow \int_{-\infty }^{\infty }r\left( \omega \right) \varphi (\omega
)d\omega \text{ \ as }N\rightarrow \infty ,
\end{equation*}
for any test function $\varphi \in D.$
\end{theorem}

\textbf{Remark.} The convergence of $r_{N}\left( \omega \right) $ to $%
r\left( \omega \right) $ follows from explicit formulas
(\ref{rnr}), but it also can be derived from a much more general
fact that the spectral measure of the Schr\"{o}dinger operator $H$
in $L^{2}(R)$ is a weak limit of the spectral measures of the
operators with truncated potentials. This fact is valid for very
general potentials without the assumption on periodicity. We shall
discuss this result elsewhere.

\textbf{Proof.} The monodromy operator $M_{\omega }$ satisfies the relation (%
\ref{s3}):
\begin{equation*}
M_{\omega }^{2}-2\cos k\left( \omega \right) M_{\omega }+I=0.
\end{equation*}
Formula (\ref{s6}) follows from here by induction. The first three relations
of (\ref{s7}) are immediate consequences of (\ref{s6}). In order to get the
fourth one we note that the eigenvalues of $M_{\omega }^{N}$ are $\mu ^{\pm
}\left( \omega \right) =e^{\pm iNk\left( \omega \right) }.$ Thus,
\begin{equation}
\alpha _{N}+\delta _{N}=\text{Tr}M_{\omega }^{N}=e^{iNk\left( \omega \right)
}+e^{-iNk\left( \omega \right) }=2\cos Nk\left( \omega \right) .  \label{s7a}
\end{equation}
\ In order to prove (\ref{ss}) we note that the relations $\det M_{\omega
}^{N}=\det M_{\omega }=1$ imply that
\begin{equation}
||M_{\omega }||^{2}-2=(\alpha -\delta )^{2}+(\beta +\gamma )^{2},  \label{m2}
\end{equation}
and a similar relation is valid for $||M_{\omega }^{N}||^{2}-2.$
From here and (\ref{s7}) it follows that
\begin{eqnarray*}
||M_{\omega }^{N}||^{2}-2 &=&(\alpha _{N}-\delta _{N})^{2}+(\beta
_{N}+\gamma _{N})^{2}=\frac{\sin ^{2}Nk\left( \omega \right) }{\sin
^{2}k\left( \omega \right) }[\left( \alpha -\delta \right) +(\beta +\gamma
)^{2}] \\
&=&\frac{\sin ^{2}Nk\left( \omega \right) }{\sin ^{2}k\left( \omega \right)
}[||M_{\omega }||^{2}-2].
\end{eqnarray*}
This justifies (\ref{ss}).

Let us prove the second statement of the Theorem. From (\ref{r1}%
) and (\ref{s7})\ it follows that
\begin{equation*}
r_{N}\left( \omega \right) =-\frac{\alpha _{N}-\delta _{N}+i\left( \beta
_{N}+\gamma _{N}\right) }{\alpha _{N}+\delta _{N}+i\left( \gamma _{N}-\beta
_{N}\right) }=-\frac{\frac{\sin Nk}{\sin k}\left( \left( \alpha -\delta
\right) +i\left( \beta +\gamma \right) \right) }{2\cos Nk+i\frac{\sin Nk}{%
\sin k}\left( \gamma -\beta \right) }.
\end{equation*}
This justifies (\ref{rnr}). In order to get (\ref{rnr1}) we note that (\ref
{ji}) implies that

\begin{equation*}
\left[
\begin{array}{c}
1+r \\
i\left[ 1-r\right]
\end{array}
\right] =c\left[
\begin{array}{c}
1 \\
m^{+}
\end{array}
\right] ,
\end{equation*}
or $\frac{1+r}{i(1-r)}=\frac{1}{m^{+}}$,\ \ and therefore,\ $\ r=\frac{%
i-m^{+}}{i+m^{+}}.$ \ \ From here and (\ref{s4a}) it follows that
\begin{eqnarray*}
r &=&\dfrac{i-\frac{e^{ik}-\alpha }{\beta }}{i+\frac{e^{ik}-\alpha }{\beta }}%
=\frac{\alpha +i\beta -e^{ik}}{e^{ik}-\left( \alpha -i\beta \right) },\text{
\ \ and} \\
r &=&\dfrac{i-\frac{\gamma }{e^{ik}-\delta }}{i+\frac{\gamma }{e^{ik}-\delta
}}=\frac{ie^{ik}-\delta i-\gamma }{ie^{ik}-\delta i+\gamma }=\frac{%
e^{ik}-\left( \delta -i\gamma \right) }{e^{ik}-\left( \delta +i\gamma
\right) }.
\end{eqnarray*}
Hence,
\begin{equation}
r=\frac{\left( \alpha +i\beta \right) -\left( \delta -i\gamma \right) }{%
2e^{ik}-\left( \alpha -i\beta \right) -\left( \delta +i\gamma \right) }=%
\frac{\left( \alpha -\delta \right) +i\left( \beta +\gamma \right) }{%
2e^{ik}-(\alpha +\delta )+i\left( \beta -\gamma \right) }=\frac{\left(
\alpha -\delta \right) +i\left( \beta +\gamma \right) }{2i\sin k+i\left(
\beta -\gamma \right) }  \label{fin}
\end{equation}
The last equality is a consequence of (\ref{cos}), and it implies
(\ref{rnr1}).

The third statement of the Theorem follows immediately from (\ref{p1}) and (%
\ref{ss}). The analyticity of $r_{N}\left( \omega \right) $ and $r$ in $%
C^{+} $ and their continuity in $\overline{C^{+}}$ follow from the explicit
formulas (\ref{rnr}), (\ref{rnr1}). Furthermore, if Im$%
k\left( \omega \right) >0, $ then
\begin{equation*}
\frac{\cos Nk(\omega )}{\sin Nk(\omega )}\rightarrow -i\text{ \ \ as }%
N\rightarrow \infty ,
\end{equation*}
and this justifies the convergence of $r_{N}\left( \omega \right) $ to $r$
when $\omega \in C^{+}.$ The weak convergence on the real axis is an obvious
consequence of the convergence in the complex half plane.

The proof of Theorem \ref{t1} is complete.

\textbf{III. Relation between slowing down and reflection.} We
start this section with considering the scattering problem for the
operator $H_{N}. $ The group velocity (for the periodic in $R$
potential $q$) is defined (see (\ref{asf})) as
\begin{equation}
V_{g}=L/k^{\prime }\left( \omega \right) .  \label{vg}
\end{equation}
In spite of the fact that formula (\ref{vg}) concerns periodic in
$R$ potentials, one can expect that a finite slab of a periodic
medium with a small group velocity, defined by (\ref{vg}), can be
used as a device for slowing down the propagation of wave packets.
However, it could happen that the wave packets with frequencies
around $\omega =\omega _{0}$, for which $V_{g}\left( \omega
_{0}\right) $ is very small, will not enter the slab, but will be
reflected almost
completely. Hence, one needs to find the situations when simultaneously $%
V_{g}$ is small and the reflection is not very large. To study
this problem we consider the solutions of the stationary
scattering problem for the operator $H_{N}$ with values of $\omega
=\omega _{0}$, for
which $V_{g}$ is small, and evaluate the reflection coefficient $%
r_{N}(\omega _{0}).$ When it is needed, we shall make more
rigorous analysis of the situation by studying the corresponding
time dependent problem.

Let us discuss the reasons for the group velocity to be small. Due to (%
\ref{cos}), the dispersion $k=k\left( \omega \right) $ on each
band $b_{n}$ is equal to $\pm $ $\arccos F\left( \omega \right) $,
where $F\left( \omega \right) =\frac{1}{2}$Tr$M_{\omega }.$ The
sign in this representation depends on the sign of $F^{\prime
}\left( \omega \right) $. Hence,
\begin{equation}
k^{\prime }\left( \omega \right) =\frac{\pm F^{\prime }\left( \omega \right)
}{\sqrt{1-F^{2}\left( \omega \right) }},  \label{kpr}
\end{equation}
and $k^{\prime }\left( \omega \right) $ can be large either because $\omega $
is close to a value $\omega _{0}$ for which $|F\left( \omega _{0}\right) |=1$
or because $\left| F^{\prime }(\omega )\right| >>1$. In the first case,
there are still two different possibilities: $F^{\prime }\left( \omega
_{0}\right) \neq 0$ or $F^{\prime }\left( \omega _{0}\right) =0.$ If $%
F^{\prime }\left( \omega _{0}\right) \neq 0$ then $\omega _{0}$ is
a non-degenerate band edge (see Lemma \ref{la}). If $F^{\prime
}\left( \omega _{0}\right) =0$ then $\omega _{0}$ is a degenerate
band edge: two bands are adjacent at $\omega _{0}$, but an
arbitrarily small perturbation of the operator may open a gap
around $\omega _{0}.$ One can consider $\omega _{0}$ as a
degenerate gap consisting of one point (see Fig. 1). Thus, we
should study the following three cases: the spectrum of the
incident wave packet is supported in a neighborhood of a
non-degenerate band edge, in a neighborhood of a degenerate band
edge or in a region where $\left| F^{\prime }(\omega )\right|
>>1$.

\textit{Wave packets with frequencies in a gap or near a
non-degenerate band edge.} Note that the transmission coefficient
$t_{N}(\omega )$ can not vanish at a fixed $\omega $ if $N<\infty
,$ since otherwise the scattering solution is zero for $x>NL,$ and
therefore has to be identically zero. A part of the energy of an
arbitrary incident wave packet will be transmitted through a
finite slab of any periodic medium. However, the following
statement shows that the transmission will be negligibly small if
the slab is long ($N$ is big) and the frequency spectrum of the
incident wave packet is supported in a gap of the corresponding
periodic problem or near a non-degenerate band edge. Let us recall
that if $\omega _{0}$ is a band edge, then $|\cos k(\omega
_{0})|=1,$ i.e. $k(\omega _{0})=n\pi ,$ $n$ is an integer.

\begin{theorem}
\label{tm}1). If $\omega $ belongs to a gap of the periodic Hamiltonian (\ref
{ha}) then
\begin{equation*}
|t_{N}(\omega )|=O(e^{-\sigma N})\text{ \ \ \ as }N\rightarrow \infty ,
\end{equation*}
where $\sigma >0$ depends on the distance from $\omega $ to the closest band.

2). If $\omega _{0}$ is a non-degenerate band edge and $\omega $
is so close to $\omega _{0}$ that $N|k(\omega )-k(\omega
_{0})|<\pi /2$ then
\begin{equation}
|t_{N}(\omega )|=O(\frac{1}{N})\text{ \ \ \ as }N\rightarrow \infty .
\label{vt}
\end{equation}
\end{theorem}

\textbf{Remark}. The assumption in the second statement means, roughly
speaking, that $|\omega -\omega _{0}|=O(1/N)$.

\textbf{Proof.} If $\omega $ belongs to a gap then $\mu ^{\pm
}=e^{\pm ik\left( \omega \right) }$\ are real, i.e. $k\left(
\omega \right) $ is purely imaginary modulus $n\pi .$ Thus,
\begin{equation}
\frac{|\sin Nk\left( \omega \right) |}{|\sin k\left( \omega \right) |}%
=O(e^{\sigma N})\ \text{\ \ as }N\rightarrow \infty .  \label{exp}
\end{equation}

If $\omega _{0}$ is a band edge, then $|\mu ^{\pm }(\omega _{0})|=1$ and $%
k\left( \omega _{0}\right) =n\pi .$ If $N|k\left( \omega \right) -n\pi |<\pi
/2$ then
\begin{equation}
\frac{|\sin Nk\left( \omega \right) |}{|\sin k\left( \omega \right) |}=O(N)\
\text{\ \ as }N\rightarrow \infty .  \label{exp1}
\end{equation}

The statements of the Theorem will follow from (\ref{exp}), (\ref{exp1}%
), and (\ref{tnt}) if we show that $||M_{\omega }||^{2}>2$ when
$\omega $ belongs to a gap or a neighborhood of a non-degenerate
band edge.

Formula (\ref{m2}) implies that for any $\omega$ we have $%
||M_{\omega }||^{2}\geq 2,$ and $||M_{\omega }||^{2}=2$ only if
$\alpha =\delta ,$ $\beta =-\gamma .$ In the latter case, the
eigenvalues $\mu
^{\pm }$\ of the matrix $M_{\omega }$ are equal to $\alpha \pm i\beta .$ If $%
\omega $ belongs to a gap, then $\mu ^{\pm }$\ have to be real and (\ref{aaa}%
) has to hold. This is impossible if $\mu ^{\pm }=\alpha \pm
i\beta .$\ Hence, $||M_{\omega }||^{2}>2$ in gaps. If $\omega
_{0}$ is a band edge, then $|\mu ^{\pm }(\omega _{0})|=1$. Thus,
$\alpha =\delta =\pm 1,$ $\beta =-\gamma =0$ at $\omega =\omega
_{0}.$ This contradicts Lemma \ref {la} in the case of a
non-degenerate band edge. Thus $||M_{\omega _{0}}||^{2}>2$ and
therefore this is also true for $\omega $ close enough to $\omega
_{0}.$

This completes the proof.

\textit{Wave packets with frequencies near a degenerate band edge.
}We shall show that the reflection coefficient $r_{N}$ is zero at
any degenerate band edge $\omega =\omega _{0},$ i.e. the medium is
transparent for the plane wave with the frequency $\omega _{0}$
(not for the wave packets). It follows immediately from
(\ref{kpr}) and Lemma \ref{la} that the group velocity is zero at
non-degenerate band edges. This is not true for degenerate band
edges. We shall show that the group velocity in this case is
''usually'' not small, but can be made small for a specific
medium. The next theorem provides formulas for $r_{N}$ and $V_{g}$ at $%
\omega =\omega _{0}$. It is followed by a discussion of how the
group velocity can be small, and then by a discussion of the
transparency of the media for wave packets.

\begin{theorem}
\label{ta}Let $\omega =\omega _{0}$ be a degenerate \textit{band
edge} (i.e. $F\left( \omega _{0}\right) =\pm 1,$ $F^{\prime
}\left( \omega _{0}\right) =0$ ). Then
\begin{equation*}
r_{N}(\omega _{0})=0,\text{ \ and }V_{g}(\omega _{0})=\frac{L}{\sqrt{%
|F^{\prime \prime }(\omega _{0})|}}.
\end{equation*}
\end{theorem}

\textbf{Proof.} If $\omega =\omega _{0}$ is a degenerate band edge
then (see Lemma \ref{la}) $M_{\omega _{0}}=\pm I$ and $F^{\prime
\prime }\left( \omega _{0}\right) \neq 0,$ where $I$ is the
identity matrix$.$ The first relation allows us to pass to the
limit in (\ref{rnr}) as $\omega
\rightarrow \omega _{0}.$ The numerator in the right hand side of (\ref{rnr}%
) vanishes as $\omega \rightarrow \omega _{0}.$ The denominator
converges to $\pm 2/N,$ since $k(\omega _{0})=n\pi .$ This proves
the first statement of Theorem \ref{ta}. Let us evaluate the group
velocity. We shall assume that $F\left( \omega \right) =1.$ The
other case ($F\left( \omega \right) =-1$) can be treated
similarly. Then we have:
\begin{equation*}
F\left( \omega \right) =1-c(\omega -\omega _{0})^{2}+O\left( |\omega -\omega
_{0}\right| ^{3})\text{ \ as }\omega \rightarrow \omega _{0},\text{ \ where }%
c>0.
\end{equation*}
This and (\ref{cos}) imply that $k\left( \omega \right) =2n\pi \pm (\omega
-\omega _{0})\sqrt{2c}+O(\left| \omega -\omega _{0}\right| ^{3/2}),$ where
the sign has to be chosen to satisfy (\ref{imk})$.$ Hence,
\begin{equation*}
k\left( \omega \right) =2n\pi +(\omega -\omega _{0})\sqrt{2c}+O(\left|
\omega -\omega _{0}\right| ^{3/2}),
\end{equation*}
and $k^{\prime }\left( \omega \right) =\sqrt{2c}+O(\left| \omega
-\omega _{0}\right| ^{1/2})=\sqrt{|F^{\prime \prime }(\omega
_{0})|}+O(\left| \omega -\omega _{0}\right| ^{1/2}).$ This
completes the proof of the Theorem.

Since we are looking for media with a small group velocity, we
need the value of $F^{\prime \prime }(\omega _{0})$ to be large.
Hence, the potential has to be large, since the elements of the
monodromy matrix and their derivatives with respect to $\omega $
can be estimated through the $L^{1}$-norm of the potential over
the period. We restrict ourselves to considering a natural class
of potentials of the form $q(x)=Av(x)$, where $A>0$ is a big
parameter and $v(x)$ is a piece-wise infinitely smooth periodic
function which is separated from zero on each interval of the
continuity. We shall show that the group velocity is not small if
$v(x)$ is positive, negative, or changes sign one time on the
period. However, a more sophisticated choice of $v$ may provide an
example of a potential for which $|F^{\prime \prime }(\omega
_{0})|\rightarrow \infty $ as $A\rightarrow \infty .$ We shall
provide one such an example, in which the intervals where $v$ is
positive or negative are separated by the intervals where $v=0.$
To make the calculation simpler we consider a non-smooth $v$ of
the form
\begin{equation}
q(x)=A\sum_{n=-\infty }^{\infty }[\delta (x-nl)-\delta (x-2nl)],\text{ \ \ }%
l=L/2,  \label{vd}
\end{equation}
where $\delta (x)$ is the Heviside delta function. Let us also
mention that the wave equation $u_{tt}+Hu=0$ is normalized in such
a way that the speed of light is equal to one. So, when we speak
about a small group velocity, we mean that it is small compared to
one.

\begin{theorem}
\label{tgr}Let the frequency $\omega =$ $\omega _{0}$ be fixed and let the
potential $q$ have the form $q(x)=Av(x)$, where $v$ is a piece wise $%
C^{\infty }$function, $v(x+L)=v(x),$ and $A\rightarrow \infty .$ Then

a) if $v\in C^{\infty }\ $on $(0,L)$ and $v(x)\geq c>0,$ then
$\omega _{0}$ belongs to a gap when $A$ is large enough;

b) if $v\in C^{\infty }$ on $(0,L),$ $v(x)\leq -c<0,$ and $\omega _{0}$ is a
degenerate band edge, then $V_{g}(\omega _{0})\rightarrow \infty $ as $%
A\rightarrow \infty $;

c) if $v\in C^{\infty }\ $on $(0,a)$ and on $(a,L),$ $0<a<L,$ and
$v(x)\geq c>0$ on one of these intervals, $v(x)\leq -c<0$ on
another one, then $\omega _{0}$ is not a degenerate band edge when
$A$ is large enough;

d) if $q$ has the form (\ref{vd}), then $\omega =\omega
_{n}=\frac{n\pi }{l},$ $n\geq 1,$ are non-degenerate band edges,
and $V_{g}(\omega _{0})=O(A^{-1})$ when $A\rightarrow \infty $ and
$A|\sin \omega l|<1.$
\end{theorem}

\textbf{Proof. }Let the assumptions of part a) hold. Then there are two solutions $%
\psi =\psi _{s},$ $s=1,2,$ of the equation
\begin{equation}
-\psi ^{\prime \prime }+(Av(x)-\omega ^{2})\psi =0  \label{aq}
\end{equation}
which have the following asymptotic behavior as $x\in \lbrack 0,L],$ $%
|\omega -\omega _{0}|<1,$ $A\rightarrow \infty $ $\ ($see for
example,
\cite{fed} or \cite{v})$:$%
\begin{eqnarray}
\psi _{1}(x)\text{ } &\sim &\text{ }e^{\int_{0}^{x}\sqrt{Av(t)}%
dt}\sum_{0}^{\infty }\alpha _{j}(x,\omega )A^{-j/2},\text{ \ \ }\alpha
_{0}=v^{-1/4}(x),  \notag \\
\psi _{2}(x)\text{ } &\sim &\text{ }e^{-\int_{0}^{x}\sqrt{Av(t)}%
dt}\sum_{0}^{\infty }\beta _{j}(x,\omega )A^{-j/2},\text{ \ \ }\beta
_{0}=v^{-1/4}(x).  \label{ww}
\end{eqnarray}
This allows us to write the asymptotic behavior of the Pr\"{u}ffer monodromy
matrix $M_{\omega }$ and its trace $2F(\omega ).$ We get that
\begin{equation}
F(\omega )=([\frac{v(0)}{v(L)}]^{1/4}+[\frac{v(L)}{v(0)}]^{1/4})\cosh
(\int_{0}^{L}\sqrt{Av(t)}dt)(1+O(A^{-1/2})).  \label{bv}
\end{equation}
Hence, $|F(\omega )|>1$ when $|\omega -\omega _{0}|<1,$ and $A$ is
large enough. Thus, $\omega _{0}$ belongs to a gap.

Let the assumptions of part b) hold. Formulas (\ref{ww}) still
provide asymptotics of two solutions of (\ref{aq}), which
oscillate as $A\rightarrow \infty ,$ since now $\sqrt{v}$ is a
purely imaginary function. This leads to the following analog of
(\ref{bv}):
\begin{equation}
F(\omega )=(|\frac{v(0)}{v(L)}|^{1/4}+|\frac{v(L)}{v(0)}|^{1/4})\cos
(\int_{0}^{L}\sqrt{A|v(t)|}dt)+O(A^{-1/2})),  \label{fom}
\end{equation}
and this formula can be differentiated with respect to $\omega .$ The main
term of the asymptotic expansion of $F(\omega )$ does not depend on $%
\omega ,$ and $F^{\prime \prime }(\omega )\rightarrow 0$ as $A\rightarrow
\infty .$ Together with Theorem \ref{ta} this proves that $V_{g}(\omega
_{0})\rightarrow \infty $ as $A\rightarrow \infty .$

Let the assumptions of part c) hold, and let $v\geq c>0$ on
$(0,a),$ $v\leq -c<0$ on $(a,L)$ (the other case can be considered
absolutely similarly. We shall use solutions $\psi _{1},\psi
_{2}$ in order to construct the transfer matrix $T^{1}$ over the interval $%
(0,a),$ and similar solutions, with the lower limits in the
integrals replaced by $a,$ in\ order to construct the transfer
matrix $T^{2}$ over the interval $(a,L)$. Then we arrive to the
following formula for the
main term $\widetilde{M}_{\omega }$ of the asymptotic expansion of $%
M_{\omega }$ as $A\rightarrow \infty :\widetilde{M}_{\omega
}=T^{2}T^{1}$, where
\begin{eqnarray*}
T^{1} &=&\left[
\begin{array}{cc}
\lbrack \frac{v(0)}{v(a_{-})}]^{1/4}\cosh Q & \frac{\omega }{\sqrt{A}}[\frac{%
1}{v(0)v(a_{-})}]^{1/4}\sinh Q \\
\frac{\sqrt{A}}{\omega }[v(0)v(a_{-})]^{1/4}\sinh Q & [\frac{v(a_{-})}{v(0)}%
]^{1/4}\cosh Q
\end{array}
\right] , \\
T^{2} &=&\left[
\begin{array}{cc}
|\frac{v(a_{+})}{v(L)}|^{1/4}\cos P & \frac{\omega }{\sqrt{A}}|\frac{1}{%
v(a_{+})v(L)}|^{1/4}\sin P \\
-\frac{\sqrt{A}}{\omega }|v(a_{+})v(L)|^{1/4}\sin P & |\frac{v(L)}{v(a_{+})}%
|^{1/4}\cos P
\end{array}
\right] .
\end{eqnarray*}
Here
\begin{equation*}
v(a_{\pm })=\lim_{x\rightarrow a\pm 0}v(x),\text{ \ }Q=\int_{0}^{a}\sqrt{%
Av(t)}dt,\text{ \ }P=\int_{a}^{L}\sqrt{A|v(t)|}dt.
\end{equation*}
If $\omega =\omega _{0}$ is a degenerate band edge, then $M_{\omega }$ and $%
\widetilde{M}_{\omega }$ have to be equal to the identity matrix (Lemma \ref
{la}). We evaluate non-diagonal elements of $\widetilde{M}_{\omega },$ using
formulas above, and equate them to zero. Then we arrive to a homogeneous
system of equations for $\sinh Q\cos P$ and $\cosh Q\sin P$ with a non-zero
determinant. Hence,
\begin{equation}
\sinh Q\cos P=\cosh Q\sin P=0.  \label{cvb}
\end{equation}
Since $Q>0,$ we have $\sinh Q\neq 0,$ $\cosh Q\neq 0,$ and
therefore, (\ref{cvb}) can not be valid, i.e. our assumption on
the existence of a degenerate band edge is wrong.

Now we are going to prove the last statement of the Theorem. Let
the potential be given by (\ref{vd}). Then the Pr\"{u}ffer
monodromy matrix $M_{\omega }$ is the product:
\begin{equation*}
M_{\omega }=T_{\omega }(l^{+},L^{-})T_{\omega }(l^{-},l^{+})T_{\omega
}(0^{+},l^{-})T_{\omega }(0^{-},0^{+}),
\end{equation*}
where$\ T_{\omega }(a,b)$ is the transfer matrix which maps the Pr\"{u}ffer
Cauchy data at $x=a$ of any solution to the equation $H\psi =\omega ^{2}\psi
$ into the Pr\"{u}ffer Cauchy data of the same solution at $x=b,$ and the
upper index $\pm $ in the arguments of $T_{\omega }$ indicates the limit
value of $T_{\omega }$ from the right or the left, respectively. Elementary
calculations give
\begin{equation}
M_{\omega }=\left[
\begin{array}{ll}
\cos \omega l & \sin \omega l \\
-\sin \omega l & \cos \omega l
\end{array}
\right] \left[
\begin{array}{ll}
1 & 0 \\
\frac{A}{\omega } & 1
\end{array}
\right] \left[
\begin{array}{ll}
\cos \omega l & \sin \omega l \\
-\sin \omega l & \cos \omega l
\end{array}
\right] \left[
\begin{array}{ll}
1 & 0 \\
-\frac{A}{\omega } & 1
\end{array}
\right] ,  \label{mmm}
\end{equation}
and
\begin{equation}
F(\omega )=\frac{1}{2}\text{Tr}M_{\omega }=\cos \omega L-\frac{A^{2}}{%
2\omega ^{2}}\sin ^{2}\omega l.  \label{fff}
\end{equation}
One can easily check that $M_{\omega }$ is the identity matrix at
$\omega =\omega _{n}=\frac{n\pi }{l},$ $n\geq 1.$ Hence, all these
points $\omega _{n}$ are degenerate band edges. If $A\rightarrow
\infty $ and $A|\sin \omega l|<1$ then (\ref{kpr}) and (\ref{fff})
imply
\begin{equation*}
|k^{\prime }(\omega )|=\frac{|L\sin \omega L+\frac{A^{2}}{\omega ^{2}}\sin
\omega l\cos \omega l-\frac{A^{2}}{4\omega ^{3}}\sin ^{2}\omega l|}{\sqrt{%
\sin ^{2}\omega L+\frac{A^{2}}{\omega ^{2}}\sin ^{2}\omega l-\frac{A^{4}}{%
4\omega ^{4}}\sin ^{4}\omega l}}=\frac{A}{\omega }\frac{|A^{-2}L+\cos \omega
l-\frac{1}{4\omega }|\sin \omega l||}{\sqrt{\omega ^{2}+1-\frac{A^{2}}{%
4\omega ^{2}}\sin ^{2}\omega l}}=O(A).
\end{equation*}
This justifies the last statement of Theorem \ref{tgr} and
completes the proof of the Theorem.

Theorem \ref{ta} together with the last statement of Theorem
\ref{tgr} can produce an impression that both the group velocity
and the reflection coefficient are small for the potential
(\ref{vd}) if the spectrum of the incident wave packet belongs to
a small enough neighborhood of a degenerate band edge $\omega
=\omega _{n}=\frac{n\pi }{l},$ $n\geq 1$. This impression is
wrong. Let us consider this situation in more detail, taking into
account that the wave packet does not have a fixed frequency
$\omega =\omega _{n},$ but the
frequencies of the waves in the wave packet belong to an interval $%
\Delta $ centered at $\omega _{n}.$ We assume that $\Delta $
belongs to the union of two bands adjacent at $\omega =\omega
_{n}.$ Since $|F(\omega )|\leq 1$ on these bands, (\ref{fff})
implies that the length of these two bands has order $O(A^{-1}).$
We assume that $\Delta $ covers only a small part of the bands,
since the wave packets with frequency spectrum strictly inside of
a band will be studied below. So, we assume that the length
$|\Delta |$ of $\Delta $ has order $o(A^{-1}).$

If the spectrum of a wave packet belongs to an interval of length $%
|\Delta |,$ then the majority of the energy of the packet at a fixed time ($%
t=0)$ is concentrated on a space interval of length $1/|\Delta |$.
The incident wave packet comes from $x=-\infty ,$ it propagates
with the speed one in the region where the potential is zero
(along the negative semiaxis, in particular), and it reaches the
point $x=0$ at time $t=0.$ After that, time $1/|\Delta |$ is
needed for the incident packet to enter the interval $[0,NL]$
where the potential is supported. We would like to have the time
of the delay of propagation of the packet through the interval
$[0,NL]$ to be at least of the same order. This leads to the
relation $\frac{NL}{V_{g}}=O(\frac{1}{|\Delta |}).$ In our case it
means that
\begin{equation}
N=O(\frac{A}{|\Delta |}).  \label{nmn}
\end{equation}
Thus, we shall assume that $|\Delta |=o(A^{-1})$ and $N$ satisfies
(\ref{nmn})$.$

\begin{theorem}
Let the potential $q$ have the form (\ref{vd}), where $%
A\rightarrow \infty ,$ let $\omega \in \Delta $ where $\Delta $ is an
interval centerd at $\omega =\omega _{n}=\frac{n\pi }{l},$ $n\geq 1,$ and
let $|\Delta |=o(A^{-1}).$ Then
\begin{equation*}
|t_{N}\left( \omega \right) |^{2}\sim \frac{4}{\omega ^{-4}A^{2}\sin
^{2}Nk\left( \omega \right) +4},
\end{equation*}
where $\cos k(\omega )=F\left( \omega \right) $ with $F\left( \omega \right)
$ given in (\ref{fff}).
\end{theorem}

\textbf{Remark.} One can see that the complete transparency ($|t_{N}\left(
\omega \right) |=1$) takes place when $\sin Nk\left( \omega \right) =0.$
There is a finite number (of order $o(N)$) of such points $\omega $ in $%
\Delta ,$ and $\omega =\omega _{n}$ is one of them$.$ Some very
small
neighborhoods of these points are also transparent enough. However, $%
|t_{N}\left( \omega \right) |\sim 1/A$ on the major part of
$\Delta$. This indicates that the majority of the energy of the
incident pulse will be reflected. Perhaps, a more rigorous
analysis of the situation is needed to justify the last statement.
We shall provide this analysis below in a similar situation of the
wave packets with the frequency spectrum in a narrow band. In this
subsection, we shall restrict ourselves to considering the case of
$N=\infty .$

\textbf{Proof. }One can easily evaluate the product (\ref{mmm}) and then
find $||M_{\omega }||^{2}-2.\ $We get
\begin{equation*}
||M_{\omega }||^{2}-2=(\frac{A^{4}}{\omega ^{4}}+\frac{4A^{2}}{\omega ^{2}}%
)\sin ^{2}\omega l.
\end{equation*}
Since $\omega \in \Delta ,$ we have that $A\sin \omega l\rightarrow 0.$
Using (\ref{fff}) we get
\begin{equation}
\sin ^{2}k\left( \omega \right) =1-F^{2}\left( \omega \right) \sim \frac{%
A^{2}}{\omega ^{2}}\sin ^{2}\omega l\text{ \ \ as }A\rightarrow \infty ,%
\text{ \ }A\sin \omega l\rightarrow 0.  \label{sink}
\end{equation}
The statement of the Theorem follows from (\ref{tnt}) and the last two
formulas.

We shall conclude this subsection by considering the limiting case
of $N=\infty $ (the potential is supported on semiaxis $x\geq 0)$,
which supports the conclusion above that the majority of the
energy of the incident wave packet is reflected.

\begin{theorem}
\label{t4}Let
\begin{equation*}
H=-\frac{d^{2}}{dx^{2}}+A\sum_{n=0}^{\infty }[\delta (x-nl)-\delta (x-2nl)],%
\text{ \ \ }A\rightarrow \infty .
\end{equation*}
Let $\Delta $ be an interval centered at $\omega _{n}=\frac{n\pi }{l}$ and $%
A|\Delta |\rightarrow 0.$ Then the following relation is valid for the
reflection coefficient:
\begin{equation*}
r(\omega )=-1+O(A^{-1}),\text{ \ \ }A\rightarrow \infty ,\text{ \ }\omega
\in \Delta .
\end{equation*}
\end{theorem}

\textbf{Proof}. Formula (\ref{rnr1}) implies
\begin{equation}
r(\omega )+1=\frac{2\sin k\left( \omega \right) +2\beta -i\left( \alpha
-\delta \right) }{2\sin k\left( \omega \right) +\left( \beta -\gamma \right)
}.  \label{r+}
\end{equation}
One can get easily from (\ref{mmm}) that
\begin{equation}
\alpha -\delta =-2\frac{A}{\omega }\sin \omega l(1+o(1)),\text{ \ }\beta =o(%
\frac{A}{\omega }\sin \omega l),\text{ \ }\beta -\gamma =\frac{A^{2}}{\omega
^{2}}\sin \omega l(1+o(1)).  \label{r-}
\end{equation}
The statement of the Theorem follows immediately from (\ref{r+}), (\ref{r-}%
), and (\ref{sink}).

\textit{Wave packets with frequencies inside of a band.} For the
same reasons as in the case of degenerate band edges, it is
natural to consider the class of potentials of the form
$q(x)=Av(x),$ where $v$ is a piece-vise smooth function,
$A\rightarrow \infty . $ If $v(x)\geq c_{0}>0$ then (see Theorem
\ref{tgr}) any fixed $\omega =\omega _{0}$ belongs to a gap when
$A$ is large enough. This case is studied above (see Theorem
\ref{tm}). If $v(x)\leq -c_{0}<0$ and $v$ is smooth enough then
$F(\omega )$ is given by (\ref{fom}). Thus, $|F(\omega )|$
can be made less than one when $A$ is large, and the chosen $%
\omega $ will belong to a band. However, (\ref{fom}) implies that $F^{\prime
}(\omega )=O(A^{-1/2})$ as $A\rightarrow \infty .$ Since $\cos k(\omega
)=F(\omega ),$ we have
\begin{equation}
k^{\prime }(\omega )\sin k(\omega )=O(A^{-1/2}),\text{ \ \ \ }A\rightarrow
\infty .  \label{kai}
\end{equation}
If $\omega $ is strictly inside of a band then $\sin k(\omega )$ is
separated from zero, and (\ref{kai}) implies that $k^{\prime }(\omega
)\rightarrow 0$ as\ $A\rightarrow \infty ,$ and therefore, the group
velocity tends to infinity.

A natural way to achieve a small group velocity is to consider a
medium with a narrow band. Since $k(\omega )$ changes by $\pi $
over a band, $k^{\prime }(\omega )$ has an order $O(\pi /d)$ in a
band of width $d$.
Hence, the group velocity will have order $O(Ld/\pi )$ and it is small when $%
d$ is small. The easiest way to get a narrow band is to consider a potential
$q(x)=Av(x),$ where $v=0$ on some interval $(a,b)\subset (0,L)$ and $%
v(x)\leq -c_{0}<0$ on the remaining part of the interval of
periodicity. Solutions $\psi \in L^2$ of the equation
\begin{equation*}
H\psi -\lambda \psi =f,\text{ \ where\ \
}H=-\frac{d^{2}}{dx^{2}}+Av(x), \text{ \ Im}\lambda \neq 0,  f \in
L^2,
\end{equation*}
tend to zero on the set $\{x:v(x)\neq 0\}$ as $A\rightarrow \infty
.$ This leads to the fact that the spectrum of the operator $H$
when $A$ is
large enough, is close to the spectrum of the operator $H_{0}=-%
\frac{d^{2}}{dx^{2}} $ on the space $L^{2}(U)$ with the Dirichlet boundary
conditions on the boundary of $U,$ where the set $U$ consists of the
interval $(a,b)$ and all its shifts by a multiple of $L$. Thus, the spectrum
of $H_{0}$ is the set of eigenvalues of the Sturm-Liouville problem on $%
(a,b) $ of infinite multiplicity. Since the spectrum of $H$ has a
band-gap structure and the spectrum is close to the spectrum of
$H_{0},$ the bands of $H$ are narrow when $A$ is large enough.

Hence, we have a class of potentials for which the operator $H$
has narrow bands, and therefore, the group velocity for
frequencies inside these narrow bands is small. It remains to find
out if long but finite slabs of such media are transparent enough
(i. e. the reflection is not too large). Unfortunately, the answer
to the question above is negative: the smaller group velocity is
in these narrow bands, the larger reflection is. This situation is
similar to what we saw for the potential (\ref{vd}) near a
degenerate band edge.

In order to make arguments simpler we restrict ourselves to
considering a model case when the Hamiltonian has the form
\begin{equation}
H=-\frac{d^{2}}{dx^{2}}+A\sum_{n}\delta (x-nL),  \label{pot}
\end{equation}
which was studied in \cite{mv}. Note, that the interval $(a,b)$ in
this case coincides with $(0,L).$ Similarly to (\ref{mmm}) we get
\begin{equation}
M_{\omega }=\left[
\begin{array}{ll}
\cos \omega L & \sin \omega L \\
-\sin \omega L & \cos \omega L
\end{array}
\right] \left[
\begin{array}{ll}
1 & 0 \\
\frac{A}{\omega } & 1
\end{array}
\right] =\left[
\begin{array}{ll}
\cos \omega L+\frac{A}{\omega }\sin \omega L & \sin \omega L \\
\frac{A}{\omega }\cos \omega L-\sin \omega L & \cos \omega L
\end{array}
\right] ,  \label{mom}
\end{equation}
and
\begin{equation}
F(\omega )=\frac{1}{2}\text{Tr}M_{\omega }=\cos \omega L+\frac{A}{2\omega }%
\sin \omega l.  \label{fmm}
\end{equation}

\begin{theorem}
\label{tt}1). If $L$ and $n$ are fixed, $A\rightarrow \infty ,$ then bands $%
b_{n}$ for the Hamiltonian (\ref{pot}) have the form:
\begin{equation*}
b_{n}=[\omega _{2n-1},\omega _{2n}]=[\frac{n\pi }{L}(1-\frac{4}{AL}%
+O(A^{-2})),\text{ }\frac{n\pi }{L}],\text{ \ \ \ }n=1,2...\text{ \ }.
\end{equation*}
2). The group velocity $V_{g}(\omega )=L/k^{\prime }(\omega )$ on each band $%
b_{n}$ has the estimate
\begin{equation}
|V_{g}(\omega )|\leq \frac{2n\pi }{AL}(1+O(A^{-1})),\text{ \ }A\rightarrow
\infty .  \label{vg1}
\end{equation}
3). The absolute value of the transmission coefficient for $\omega \in b_{n}$
is equal to
\begin{equation*}
\left| t_{N}(\omega )\right| =\frac{1}{\sqrt{1+|A\frac{\sin Nk(\omega )}{%
2\omega \sin k(\omega )}|^{2}}},
\end{equation*}
and it has the order $O(A^{-1})$ for all $\omega \in b_{n},$
$n\geq 1,$ except for very narrow (of order $O(A^{-2})$)
neighborhoods of $N-1$ transparency points where
\begin{equation*}
\sin Nk(\omega )=0,\text{ \ \ }\sin k(\omega )\neq 0.
\end{equation*}
\end{theorem}

\textbf{Proof.} The inequality $|F(\omega )|\leq 1$ can be valid only if $%
\omega L$ is close to $n\pi ,$ $n\geq 1$ (since $\omega ^{-1}\sin
\omega L$ has to be small). In order to find the bands, we put
$\omega L=n\pi +\sigma , $ $|\sigma |<<1.$ Then $\cos \omega
L=(-1)^{n}\cos \sigma ,$ $\sin \omega L=(-1)^{n}\sin \sigma $, and
the inequality $|F|\leq 1$ takes the form
\begin{equation}
|g(\sigma )|\leq 1,\text{ \ \ }|\sigma |<<1,\text{\ \ \ \ \ where }g(\sigma
)=\cos \sigma +\frac{AL}{2(n\pi +\sigma )}\sin \sigma .  \label{in}
\end{equation}

Note that $g(0)=1$ and $g^{\prime }(\sigma )>c_{0}A>0$ when $|\sigma |<\pi
/4 $ and $A\rightarrow \infty .$ From here it follows that (\ref{in}) holds
if $-\varepsilon >\sigma >0$ where $\varepsilon $ is the solution of the
equation $g(\varepsilon )=-1,$ i.e.
\begin{eqnarray*}
0 &=&1+\cos \varepsilon -\frac{AL}{2(n\pi -\varepsilon )}\sin \varepsilon
=2\cos ^{2}\varepsilon /2-\frac{AL\sin \varepsilon /2\cos \varepsilon /2}{%
n\pi -\varepsilon } \\
&=&\cos \varepsilon /2[2\cos \varepsilon /2-\frac{AL\sin \varepsilon /2}{%
n\pi -\varepsilon }].
\end{eqnarray*}
Thus $\varepsilon =\varepsilon _{n}$ is the solution of the equation
\begin{equation}
\tan \varepsilon /2=\frac{2(n\pi -\varepsilon )}{AL},  \label{ep}
\end{equation}
i.e. $\varepsilon _{n}=\frac{4n\pi }{AL}+O(A^{-2}).$ Since $b_{n}=[(n\pi
-\varepsilon _{n})/L,$ $n\pi /L],$ the first statement of Theorem \ref{tt}
is proved.

Let us prove the second statement. One can easily check that
\begin{equation}
F^{\prime }(\omega )=(-1)^{n}\frac{AL^{2}}{2n\pi }(1+O(A^{-1}))\text{ \ when
}\omega \in b_{n},\text{ \ }A\rightarrow \infty .  \label{gr1}
\end{equation}
Then from (\ref{cos}) it follows that
\begin{equation}
k^{\prime }(\omega )\sin k(\omega )=(-1)^{n+1}\frac{AL^{2}}{2n\pi }%
(1+O(A^{-1}))\text{ \ when }\omega \in b_{n},\text{ \ }A\rightarrow \infty ,
\label{gr2}
\end{equation}
which implies the second statement of the Theorem. The last statement
follows immediately from (\ref{tnt}) and (\ref{mmm}).

Theorem \ref{tt} is proved.

Theorem \ref{tt} implies that if the frequency spectrum of an
incident pulse belongs to a band $b_{n}$ of the operator
(\ref{pot}) then only a small part (of order $O(A^{-2})$) of the
energy of the pulse will propagate through the medium, and the
main part will be reflected. This statement needs to be justified
more rigorously, in particular, because of the existence of the
transparency zones. Besides, the question about the speed of
propagation of the transmitted part of the energy remains to be
addressed, since the usual arguments about the group velocity
concern the propagation of wave packets in infinite media, not
through a finite slab of a periodic medium. The next theorem deals
with both problems. In particular, we shall show that the main
part of the incident pulse will be reflected. Only a part of the
energy of the pulse of order $O(A^{-2})$ enters the media. It
propagates with the group velocity, and only a negligible part of
order $O(A^{-\infty }) $ moves faster.

The following non-stationary problem describes the propagation of the pulse $%
u=u_{0}(x-t)$ through a finite slab of the periodic medium:
\begin{equation}
u_{tt}=u_{xx}-A\sum\limits_{n=0}^{N-1}\delta \left( x-nL\right) u,\text{ \ \
}(t,x)\in R^{2};\text{ \ \ \ }u=u_{0}(x-t),\text{ \ \ }t\leq 0,\text{ \ }
\label{np}
\end{equation}
where
\begin{equation}
u_{0}(x)=B^{-1/2}\alpha (x/B)e^{ix\omega _{0}},\text{ \ }\alpha \in
C_{0}^{\infty }([-1,0]),\text{ \ }B=O(A^{1+\varepsilon }),\text{ \ \ }%
\varepsilon >0,  \label{20}
\end{equation}
and $\omega _{0}$ belongs to $b_{n}.$ \ To be more exact, we assume that
\begin{equation}
\omega _{0}=\frac{n\pi }{L}(1-\frac{\theta }{AL}),\text{ \ \ \ }0<\theta <4%
\text{, \ \ }n\geq 1,  \label{21}
\end{equation}
where $\theta $ does not depend on $A$. The form of $u_{0}$ and the meaning
of the restrictions (\ref{20}), (\ref{21}) will be explained below, but
first let us mention that the energy conservation law holds for the
solutions $u$ of (\ref{np}): the energy
\begin{equation*}
E_{u}=\int_{-\infty }^{\infty
}(|u_{x}|^{2}+|u_{t}|^{2})dx+A\sum_{m=0}^{N-1}|u(mL,t)|^{2}
\end{equation*}
does not depend on time. Since the energy of the solution at $t=0$ is equal
to
\begin{equation}
E_{u_{0}}=2\int_{-\infty }^{0}|\frac{d}{dx}u_{0}{}(x)|^{2}dx,  \label{u1}
\end{equation}
we have
\begin{equation*} E_{u}=E_{u_{0}}
\end{equation*}
for all $t\geq 0.$

The function $u_{0}$ is normalized in such a way that the energy of the
incident wave $u_{0}(x-t)$ at $t=0$ is separated from zero and infinity and
tends to a constant as $B\rightarrow \infty :$
\begin{equation}
E_{u_{0}}=2\omega _{0}^{2}||\alpha (x)||_{L^{2}}^{2}+O(B^{-1}).  \label{u0}
\end{equation}
The Fourier transform of $u_{0}$ is equal to
\begin{equation}
\widetilde{u}_{0}(\omega )=B^{1/2}\widetilde{\alpha }(B(\omega -\omega
_{0})).  \label{uv}
\end{equation}
From Theorem \ref{tt} and (\ref{21}) it follows that $\omega _{0}$
belongs to $b_{n}$ and is located not very close to the end points
of $b_{n}.$ The function $\widetilde{u}_{0}$ is supported, mostly,
on a small interval\
\begin{equation}
\Delta =\{\omega :\text{ }|\omega -\omega _{0}|\leq 1/A^{1+\varepsilon /2}\}.
\label{r3}
\end{equation}
located inside of $b_{n}.$\ In fact, from Theorem \ref{tt} and (\ref{21}) it
follows that
\begin{equation}
\Delta \subset b_{n}\text{ \ and \ }\rho (\Delta ,\partial b_{n})>\beta /2A%
\text{ \ \ if }A>(2L/\beta )^{2/\varepsilon },  \label{r4}
\end{equation}
where $\beta =\min (\theta ,4-\theta )$ and $\rho (\Delta
,\partial b_{n})$ is the distance between $\Delta $ and the end
points of $b_{n}$. Furthermore, since $\alpha \in C_{0}^{\infty
},$ we have $\widetilde{\alpha }(\omega )=O(|\omega |^{-\infty })$
as $|\omega |\rightarrow \infty .$ \ This and the assumption
$B=O(A^{1+\varepsilon })$
imply the following estimate, which holds when $\omega \notin \Delta :$%
\begin{equation}
|\widetilde{u}_{0}(\omega )|=|B^{1/2}\widetilde{\alpha }(B(\omega -\omega
_{0}))|=B^{1/2}O((\frac{A^{1+\varepsilon /2}}{B})^{\infty })=O(A^{-\infty }),%
\text{ \ \ }A\rightarrow \infty .  \label{r7}
\end{equation}
Hence, $u_{0}$ can be considered as a wave packet with frequencies from $%
b_{n}$. Let us mention that it is impossible to choose such a $u_{0}$ that $%
\widetilde{u}_{0}$ is completely supported on $b_{n},$ since $u_{0}$ would
be analytic in $x$ in this case and could not be equal to zero for $x\geq 0,$
while the incident wave $u_{0}(x-t)$ is supported on the semi-axis $x<0$
when $t<0.$

We shall assume that $N$ is not very small, since otherwise it
would be unnatural to discuss a delay in the propagation of the
incident pulse of space length $B=O(A^{1+\varepsilon })\rightarrow
\infty $. One could consider values of $N$ which are independent
of $A$ and $B$, but we shall assume that $N=O(A^{\gamma })$ with
some $\gamma >\varepsilon .$ Then the final estimates will not
depend on $N,$ so the result will be simpler.

Let us evaluate the group velocity $V_{g}(\omega )=L/k^{\prime }(\omega )$
for wave packets with frequencies from $\Delta .$ From (\ref{fmm}), (\ref{21}%
) and (\ref{r3}) it follows that
\begin{equation*}
\cos k(\omega )=(-1)^{n}(1-\frac{\theta }{2})+O(A^{-\varepsilon /2})\text{ \
\ when }\omega \in \Delta .
\end{equation*}
Hence,
\begin{equation}
|\sin k(\omega )|=\frac{1}{2}\sqrt{4\theta -\theta ^{2}}+O(A^{-\varepsilon
/2})\text{ \ \ when }\omega \in \Delta .  \label{sin}
\end{equation}
From here and (\ref{gr2}) it follows that
\begin{equation*}
|V_{g}(\omega _{0})|=\frac{n\pi \sqrt{4\theta -\theta ^{2}}}{AL}%
(1+O(A^{-\varepsilon /2})).
\end{equation*}
Thus, the time
\begin{equation*}
t=t_{0}(1+O(A^{-\varepsilon /2})), \text{ \ \ where }
t_{0}=\frac{ANL^{2}}{n\pi \sqrt{4\theta -\theta ^{2}}},
\end{equation*}
is needed for a signal to pass through the interval $[0,NL]$ of the periodic
medium, if the signal propagates with the group velocity $V_{g}(\omega
_{0}). $

\begin{theorem}
(see \cite{mv}) The following estimates are valid for the solution $u$ to
the problem (\ref{np}) when $t<t_{0}(1-\delta ),$ $\delta >0,$ and $%
A\rightarrow \infty $:

1) $\ u=u_{0}(x-t)-u_{0}(-x-t)+g(x+t)$ \ \ when $x<0,$ \ where \
\begin{equation*}
0<C_{1}A^{-2}<\int_{-\infty }^{\infty }|g^{\prime
}(x)|^{2}dx<C_{2}A^{-2}<\infty .
\end{equation*}

2)
\begin{equation*}
\overset{\infty }{\underset{NL}{\int }}(u_{x}^{2}+u_{t}^{2})dx=O(A^{-\infty
})
\end{equation*}
\end{theorem}

\end{document}